\documentclass[12pt]{article}
\usepackage{CJKutf8}
\usepackage[letterpaper,top=1in, bottom=1in, left=1in, right=1in]{geometry}
\usepackage{amsmath}
\usepackage{amsfonts}
\usepackage{amssymb}
\usepackage{hyperref}
\usepackage{tikz}
\newtheorem{prop}{Proposition}
\newtheorem{conj}{Conjecture}
\newcommand{\NEpath}[4]{
	\fill[white!25]  (#1) rectangle +(#2,#3);
	\fill[fill=white]
	(#1)
	\foreach \dir in {#4}{
		\ifnum\dir=0
		-- ++(1,0)
		\else
		-- ++(0,1)
		\fi
	} |- (#1);
	\draw[help lines] (#1) grid +(#2,#3);
	\draw[dashed] (#1) -- +(#3,#3);
	\coordinate (prev) at (#1);
	\foreach \dir in {#4}{
		\ifnum\dir=0
		\coordinate (dep) at (1,0);
		\else
		\coordinate (dep) at (0,1);
		\fi
		\draw[line width=2pt,-stealth] (prev) -- ++(dep) coordinate (prev);
	};
}

\newcommand{\NEboundpath}[4]{
	\fill[white!25]  (#1) rectangle +(#2,#3);
	\fill[fill=white]
	(#1)
	\foreach \dir in {#4}{
		\ifnum\dir=0
		-- ++(1,0)
		\else
		-- ++(0,1)
		\fi
	} |- (#1);
	\draw[help lines] (#1) grid +(#2,#3);
	\draw[dashed] (#1) -- +(#3,#3);
		\draw[dashed] (0,3) -- +(#3-3,#3-3);
	\coordinate (prev) at (#1);
	\foreach \dir in {#4}{
		\ifnum\dir=0
		\coordinate (dep) at (1,0);
		\else
		\coordinate (dep) at (0,1);
		\fi
		\draw[line width=2pt,-stealth] (prev) -- ++(dep) coordinate (prev);
	};
}

\newcommand{\NEaboveboundpath}[4]{
	\fill[white!25]  (#1) rectangle +(#2,#3);
	\fill[fill=white]
	(#1)
	\foreach \dir in {#4}{
		\ifnum\dir=0
		-- ++(1,0)
		\else
		-- ++(0,1)
		\fi
	} |- (#1);
	\draw[help lines] (#1) grid +(#2,#3);
	\draw[dashed] (#1) -- +(#3,#3);
	\draw[dashed] (0,3) -- +(#3-3,#3-3);
		\draw[dashed] (0,4) -- +(#3-4,#3-4);
	\coordinate (prev) at (#1);
	\foreach \dir in {#4}{
		\ifnum\dir=0
		\coordinate (dep) at (1,0);
		\else
		\coordinate (dep) at (0,1);
		\fi
		\draw[line width=2pt,-stealth] (prev) -- ++(dep) coordinate (prev);
	};
}

\DeclareMathOperator{\im}{Im}
\DeclareMathOperator{\spp}{Span}
\DeclareMathOperator{\qdim}{qdim}

\author{Keke Zhang}
\title{Generalized Catalan numbers}
\date{}

\begin{document}

\maketitle
	\begin{abstract}
	A finitization of the Catalan numbers $ C_n $ can be defined as Euler characteristics of an algebraic structure. We conjecture the existence of a $ q $-deformed version of such structure, and provide evidence for the first two non-trivial cases. This $ q $-deformed version will likely categorify series arise from identities of the Roger-Ramanujan type.
	\end{abstract}
	
	\section{Catalan numbers}
	
	\subsection{Dyck Paths}

We begin by considering the lattice path on a square--- a sequence of North $ N(0, 1)  $  and East $ E(1, 0) $ steps in the first quadrant of the
$ xy $-plane, starting at the origin $ (0, 0) $ and ending at $ (n,n) $. The total number of paths is then \[ {2n\choose n}. \]

A Dyck, sometimes called a Catalan path, is a lattice path that lies above the diagonal line $ x=y $.
For example, we have the following:
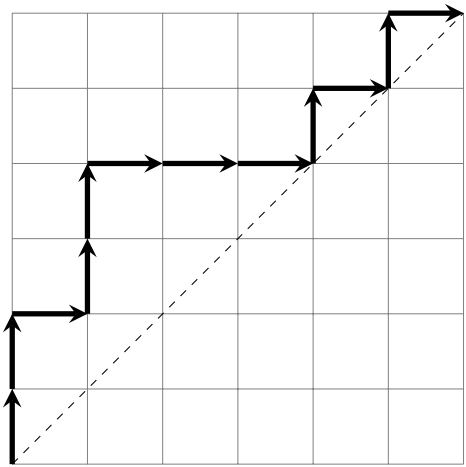
\[	\begin{tikzpicture}
\NEpath{0,0}{6}{6}{1,1,0,1,1,0,0,0,1,0,1,0};
\end{tikzpicture}
\]
Let $ C_n $ denote the nth Catalan number. We can compute $ C_n $ by substracting all pathes going across the diagonal, which, by a reflection of paths after the first one that passed the diagonal line. For example, in the picture below, the path on the left is reflected to the path on the right.

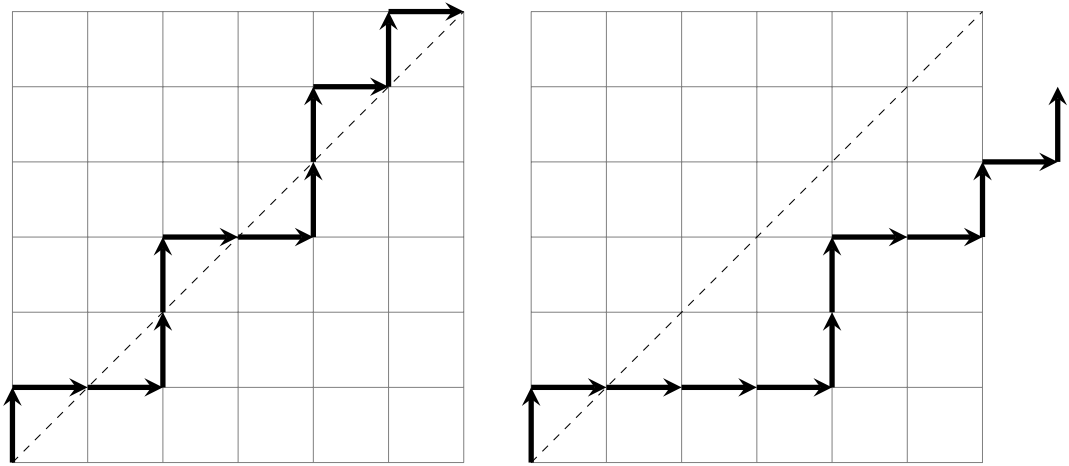
\[\begin{tikzpicture}
\NEpath{0,0}{6}{6}{1,0,0,1,1,0,0,1,1,0,1,0};
\end{tikzpicture}
\quad \quad
\begin{tikzpicture}
\NEpath{0,0}{6}{6}{1,0,0,0,0,1,1,0,0,1,0,1};
\end{tikzpicture} \]
In fact, we know that there is a one-to-one correspondence of the two types of pathes, hence, we are substracting
\[{2n\choose n-1}.\]
Then we obtain \[ C_n = {2n \choose n} -{2n\choose n-1}. \]
	
	\subsection{The Pascal Triangle}
	
We observe that the Catalan numbers are the differences of the terms in the middle column and the one next to it in a Pascal triangle, in the even rows.
	
	\[\begin{tikzpicture}
	\foreach \n in {0,...,4} {
		\foreach \k in {0,...,\n} {
			\node at (\k-\n/2,-\n) {${\n \choose \k}$};
		}
	}
	\end{tikzpicture}\]
	Where we have
	\begin{equation} \label{induction}
		 {n\choose k}={n-1\choose k-1} +{n-1\choose k} .
	\end{equation}

	Note that we start our counting at the $ 0 $-th row.
	One way to interprete this  is that in each even row, Catalan number is the Euler characteristic of a complex \[0\rightarrow C_0 \rightarrow C_1\rightarrow 0\]
	with $ C_0 $ of dimension $ {2n \choose n} $, and $ C_1 $ of dimension $ {2n \choose n+1} $.
	We will see how this complex is constructed explicitly later.

	\section{Finitization of Catalan}
	\subsection{Paths with Two Bounds}
	We can draw a line above and parallel to the diagoal line. Then we form another bound of the paths (see below). 
	
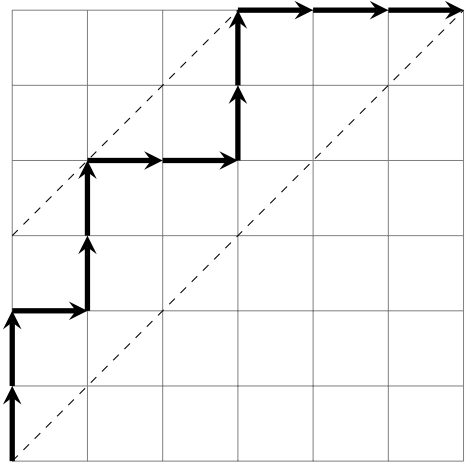
\[	\begin{tikzpicture}
\NEboundpath{0,0}{6}{6}{1,1,0,1,1,0,0,1,1,0,0,0};
\end{tikzpicture}
\]

	Then, how many paths are here?
	To answer this question, we first recall how we computed the usual Dyck paths: we take the number of all lattice paths going from $ (0,0) $ to $ (n,n) $, then substract everything goes across the diagonal line.\\
	We can count these paths with restricted height in the same manner. We draw the line right above (below) our upper (lower) bound, then, for a path goes above (below) the bound, it necessarily touches this line. We reflect the end part of this path agaist this line, then, we have shifted a path that touches the above (below) line once to a path ending at a different point.
	For instance, for a path on the left side, we reflect it to the right picture.
	\[	\begin{tikzpicture}
	\NEaboveboundpath{0,0}{6}{6}{1,1,0,1,1,0,1,1,0,0,0,0};
	\end{tikzpicture}
	\quad \quad
	\begin{tikzpicture}
	\NEaboveboundpath{0,0}{6}{6}{1,1,0,1,1,0,1,1,1,1,1,1};
	\end{tikzpicture}.
	\]
	
	In this particular example above, $ m=3 $, we have the entire path shifted to a path with endpoint $ (2,10) $. In  the general  case,  we would  have, for a path
			passes the upper bound first, 
			\[|B_{2j+1}|={2n\choose n- j(m+2)-m-1}\]
	\[|B_{2j}|={2n\choose n- j(m+2)}.\]

	Since a path can pass theupper bound  then the lower one (then the upper one again etc.).
	The ones touches the lower bound first can be computed similarly:

\[|A_{2j+1}|={2n\choose n+ j(m+2)+1}\]
\[|A_{2j}|={2n\choose n+ j(m+2)}.\]

Then, we have the total number of bounded paths to be an alternating sum 
\[\sum_{i \in \mathbb{N}}(-1)^i (|A_i|+|B_i|).\]

	One natual question to ask is that "why we are choosing the above signs when $ {2n \choose n+k}= {2n\choose n-k} $?" This question will be clear if we move to our picture of the Pascal triangle.

	For example,  for $ m=4 $, we have \[|A_1|= {2n \choose n+1} \quad |B_1|={2n \choose n-5}\]
	 \[|A_2|= {2n \choose n+6} \quad |B_2|={2n \choose n-6}\]
	 \[|A_3|= {2n \choose n+7} \quad |B_3|={2n \choose n-11}\]
	  \[|A_4|= {2n \choose n+12} \quad |B_4|={2n \choose n-12}\]
	  \[\cdots\quad\cdots\]

	  Then, on an even row of Pascal triangle,  we get an alternating sum,
	  \begin{align*}
	  \cdots+{2n \choose n-12}-{2n \choose n-11}+{2n \choose n-6}-&{2n \choose n-5}+\\
	 &  + {2n \choose n}-\\
	 &- {2n \choose n+1}  +{2n \choose n+6} - {2n \choose n+7}+{2n \choose n+12 }-	 \cdots
	  \end{align*}
	Note that with our sign choice, all $ A_i $'s are on the right side and all $ B_i $'s are on the left side. [This particular order has nothing to do with the fact that I read too many East Asian texts recently and get too used to startinging things from  right to left.]

	\subsection{An Algebraic Complex} \label{alg comp}
		The previous alternating sum naturally makes one wonder: can we construct some complexes such that this alternating sum corresponds to the Euler  characteristics of them?
	
	To answer this question, we consider the following aglebraic structure: we let $ (x_1,\cdots,x_m) $ be the generators, such that
	\[x_i^2=0\quad \forall i,\]
	\[x_ix_j =q x_j x_i \quad \forall i>j,  \] and we let
	\[q^N=1.\]
	Then, by  an expansion of the product, we have 
	\[(x_1+\cdots+x_m)^N=0,\]
	because the $ q $-binomial term includes $ (1-q^N) $ in the numerator.
	
We consider the complexies to be \begin{align*}
C;&\\
\spp \{x_1, \cdots,x_m \};&\\
\spp \{x_1x_2, \cdots, x_{m-1}x_m \};&\\
\cdots&\\
Cx_1x_2\cdots x_m.&
\end{align*}
That is, each is of dimension ${ m \choose k }$.

	We define $ \sigma  $ to be the left multiplication by \[\sum_{i=1}^{m}x_i \]

	Now we can define an algebraic complex  by taking the differentials to be \[ \sigma, \sigma^{N-1} ,\sigma, \sigma^{N-1} ,\cdots  ,\]
	such that the complex includes the term with dimension  $ {2n \choose n} $ and $ {2n \choose n+1} $.
	That is, we start from the middle, so that the following appears in the complex:
	\[ C_{[0]} \xrightarrow{ d_{[0]}} C_{[1] }\]
	where
	\[C_{[0] }= \spp \{x_{i_1}x_{i_2}\cdots x_{i_n}  \} \]
		\[C_{[1]} = \spp \{x_{i_1}x_{i_2}\cdots x_{i_{n+1}}  \} \]
		\[d_{[0]}=\sigma.\]
	Note that we picked the differentials so that $ d_{[i]}= \sigma $ for $ i $ even and $ d_{[i+1]} = \sigma^{N-1}$, then we necessarily have that $ d^2=0 $.
	We defined the indecides so that it makes more sense in our complex, and allow negative indecies to appear.
	
	In fact, we can compute to show that for everything other than the top degree, we have $ \ker d= \im d $, hence we can show that all homologies lie in one degree.
	
	This complex has the Euler characteristic of the alternating sum that we were looking at. For instance, the explicit example in the previous subsection can be obtained by setting $ N=6 $.
	
	Now, if we let $ N\mapsto \infty $, all the other components of the differential vanish, except for $ d_{[0]} $. Hence, we obtain a complex with homology equal to the Catalan numbers.

\subsection{Shifting the Lower Bound}	
	Note that  we can in fact shift this alternating sum so that it does not have to pass the middle line. This corresponds to numbers of Dyck paths with two bounds and the lower bound can be shifted to a general position below the diagonal $ x=y $. And we can get an alternating sum
	
	This can be still obtained by path reflection and with modifications:
	\[\sum_{i \in \mathbb{N}}(-1)^i (|A_i|+|B_i|),\]
	where
	\[|A_{2j}|=|B_{2j}|={2n\choose n+ j(m+s+2)}.\]
		\[|A_{2j+1}|={2n\choose n+ j(m+s+2)+s+1}\]
		\[|B_{2j+1}|={2n\choose n- j(m+s+2)-m-1},\]
where $ s $ is the distance that the lower bound shifted away from the diagonal.
	And, in fact the algebraic complex can be define  in this case similarly, by taking the alternating sums of a group of different  terms in any row of the Pascal triangle, with the differentials taken as
	\[\sigma^{m+1},\sigma^{s+1},\sigma^{m+1},\sigma^{s+1},\cdots\]
	and setting $ N=m+s+2 $, so that \[\sigma^{m+s+2}=0.\]
	Hence, the previous case can be seens as when we take $ s=0 $.
	By induction on the Pascal's triangle with relation \eqref{induction}, in fact we can see some interesting results here:
	\begin{prop} \label{prop1}
	For $ N=2 $, the alternating sum constantly vanish;\\
		for $ N=3 $,  the alternating sum gives $ 0, 1, -1 $;\\
			for $ N=4$,  the alternating sum gives $ \pm 2^n $;\\
				for $ N=5 $,  the alternating sum gives $ \pm $Fibonacci numbers;\\
					for $ N=6 $,  the alternating sum gives $ \pm 3^n, \pm\lfloor  \frac{3^n}{2} \rfloor,  \pm\lceil  \frac{3^n}{2}\rceil $ \ldots\\
					And all the above is true for different partitions into two parts of $ N $.
	\end{prop}

	\section{First Generalization: $ 3 $-Pascal's Triangle}
		\subsection{Irreducible Representations  of $ \mathfrak{sl}_2 $}
	Note that there is a deficiency in our previous example. Although our complex can be extended to the entire Pascal's triangle, Catalan numbers are only defined on the even rows. Hence, a natural question arises: how can we find a generalization of this Pascal's triangle so that we can understand the odd rows as well?
	
	We recall that there is exactly one $ n $-dimensional irreducible representation of $ \mathfrak{sl}_2 $ for each dimension $ n $, with character \[q^{-(n-1)}+ q^{-(n-3)}+\cdots+q^{(n-1)}.\] 
In particular, for the two dimensional representation, we have the character as \[q^{-1}+q^1.\]
	We consider the polynomial \[(q^{-1}+q^1)^n, \] then the coeffiencients in front of a term $ q^k $ looks like $ n \choose |k| $, which is an element in the usual Pascal's triangle. In fact, the difference between the coefficients of the middle two terms $ {2n\choose n} $ and $ {2n \choose n+1} $ corresponds to the number of $ 1 $-dimensional components in the $ 2n $-tensor of the $ 2 $-dimensional irreducible representation of $ \mathfrak{sl}_2. $

Following this intuition, we can see that our natural generalization would be using the three dimensional representation of $ \mathfrak{sl}_2 $ to construct a $ 3 $-Pascal's triangle.
We consider the polynomial obtained from multiplying the character corresponding  to it \[(q^{-2}+q^0+q^2)^n.\]
Note that the coeffiecient in front of any term $ q^{2k} $ would look like
	\[\sum_i {n \choose i, i+k } ,\] where 
	\[ {n\choose i,i+k }= \frac{n!}{i!(i+k)! (n-2i-k)!}. \]

One can see  that this corresponds to terms in a Pascal's-triangle-like structure where each term is the sum of the three terms right above it.

\[ \begin{matrix}
&	 & & & 1 &&&&\\
&& &	1 & 1 & 1&&&\\
&&1&2& 3  &2&1&&\\
&1&3&6&7&6&3&1&\\
1&4&10&16&19&16&10&4&1\\
&&&&\cdots&&&&
\end{matrix} \]
	
	\subsection{Generalized Path counting}
	
	Similar to the ususal Dyck paths, we consider the generalized lattice paths genrated by the three moves, from the origin to \[(0,2), (1,1),(2,0),\]
	respectively. We call them the$ (0,2), (1,1),(2,0) $ paths respectively.
	
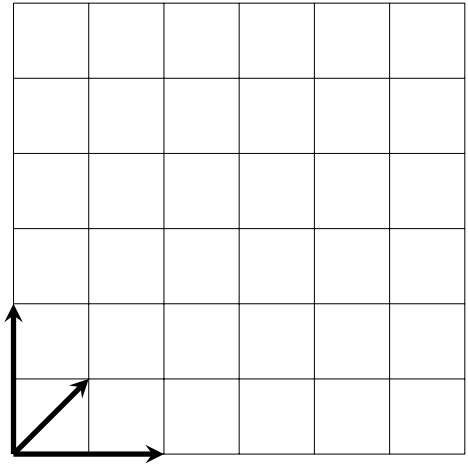
\[ \begin{tikzpicture}[domain=0:6] 
	\draw[very thin,] (0,0) grid (6,6);
	\draw[line width=2pt,-stealth] (0,0) -- (2,0) ; 
	\draw[line width=2pt,-stealth] (0,0) -- (1,1) ; 
	\draw[line width=2pt,-stealth] (0,0) -- (0,2) ; 
	\end{tikzpicture}\]

	Note that we can generalize the restricted path-counting, and in  particular, Catalan numbers, by  considering the  paths strictly above the diagonal. 
	That is, it can not have any components lying  on the diagonal. The first one has to start to $ (0,2) $, then we would have three choices. [In fact, if one draws this generalized lattice completely, then rotate it clockwise by $ \pi/4 $, one would see a sketch of I. M. Pei's pyramid entrance of the Louvre. But don't tell the French people, I think his design of the Suzhou Museum is way better.]
	
	By analogy, we can define the generalized Dyck paths with two bounds. For example, the one shown below on the right.
	
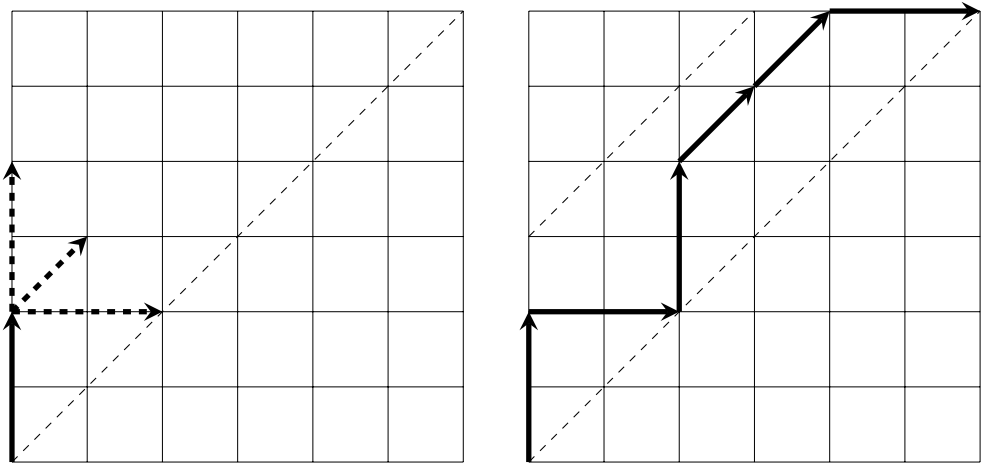
\[	  \begin{tikzpicture}[domain=0:6] 
	\draw[very thin] (0,0) grid (6,6);
	\draw[line width=2pt,-stealth] (0,0) -- (0,2) ; 
		\draw[dashed,line width=2pt,-stealth] (0,2) -- (2,2) ; 
			\draw[dashed,line width=2pt,-stealth] (0,2) -- (1,3) ; 
				\draw[ dashed,line width=2pt,-stealth] (0,2) -- (0,4) ; 
	\draw[dashed] (0,0) -- +(6,6);
	\end{tikzpicture} \quad \quad	\begin{tikzpicture}[domain=0:6] 
	\draw[very thin,] (0,0) grid (6,6);
	\draw[line width=2pt,-stealth] (0,0) -- (0,2) ; 
	\draw[line width=2pt,-stealth] (0,2) -- (2,2) ; 
	\draw[line width=2pt,-stealth] (2,2) -- (2,4) ; 
	\draw[line width=2pt,-stealth] (2,4) -- (3,5) ; 
	\draw[line width=2pt,-stealth] (3,5) -- (4,6) ; 
	\draw[line width=2pt,-stealth] (4,6) -- (6,6) ; 
	\draw[dashed] (0,0) -- +(6,6);
	\draw[dashed] (0,3) -- +(3,3);
	\end{tikzpicture}\]
	
	The number of bounded paths again correspond to the differences of the middle two columns of the $ 3$-Pascal's triangle.
	
	\[\sum_k{n\choose k,k}-\sum_k{n\choose k,k+1}. \]
	
	And similarly, we have for the double bounded paths
	\[\sum_{i \in \mathbb{N}}(-1)^i (|A_i|+|B_i|), \]
	where
	\[|A_{2j}|=|B_{2j}|=\sum_k{n\choose k,k+j(m+2) } \]
	\[|A_{2j+1}|=\sum_k{n\choose k,k+j(m+2)+1} \]
	\[|B_{2j+1}|=\sum_k{n\choose k,k-j(m+2)-(m+1)}. \]
	
	Note that both the even and the odd values of $ m $ make sense, when we  take into account of the fact that an odd value $ 2j+1 $ would include the $ (1,1) $ paths right below it, but $ 2j $ would exclude them.

Finally we remark that the above can be generalized to an $ n $-Pascal's triangle corresponding to the $ n $-dimensional irreducible representation of $ \mathfrak{sl}_2 $.

Counterintuitively, the value of the alternating sum does depend on the dimension of the representation we choose. For example, compared to proposition  \ref{prop1}, we have the following:
	\begin{prop} \label{prop2}
		In the $ 3 $-Pascal's triangle,
	For $ N=3 $, the alternating sum constantly vanish;\\
	for $ N=2 , N=4$,  the alternating sum gives $ 0, 1, -1 $;\\
	for $ N=5$,  the alternating sum gives $ \pm $Fibonacci numbers;\\
	for $ N=6$,  the alternating sum gives $ \pm $Jacobsthal numbers, that is, \[ a(n) = a(n-1) + 2a(n-2),\]  with $ a(0) = 0, a(1) = 1 .$\\
In the $ 4 $-Pascal's triangle,
	For $ N=2 ,N=4$, the alternating sum constantly vanish;\\
for $ N=3, N=5 $ ,  the alternating sum gives $ 0, 1, -1 $.\\

\end{prop}

	\section{Second Genrealization: $ q $-Pascal Triangle}
	\subsection{The $ q $-Analogs}
	We consider the $ q $-binomials
	\[{n\choose k}_q= \frac{(n!)_q}{((n+k)!)_q (k!)_q},  \]
	where $ (n!)_q $ is a factorial over the numbers \[\frac{1-q^n}{1-q}=1+q+\cdots+q^{n-1}. \]

Then, we can naturally arrange them in the following way to define the $ q $-Pascal's  triangle.
\[\begin{tikzpicture}
\foreach \n in {0,...,4} {
	\foreach \k in {0,...,\n} {
		\node at (\k-\n/2,-\n) {${\n \choose \k}_q$};
	}
}
\end{tikzpicture}\]

We note that the relation \[{n\choose k}_q={n-1\choose k-1}_q +q^{k}{n-1\choose k}_q \]
allows us to compute it inductively.

We can then define the $ q $-Catalan numbers by
\[ C_q(n)={n\choose k}_q-q{n\choose k+1}. \]

In section \ref{alg comp}, we have defined the algebraic complex, and the $ q $-binomials correspond to the  graded dimensions of each of the complexes, defined such that for each $ x_{i_1}x_{i_2}\cdots x_{i_k} $, with the indices in  ascending order, we have it in the $ q $-dimension as
\[\prod_j  q^{i_j-j} .\]
Then,  the  $ q $-dimension for the complex would be
\[\qdim=\sum_{i}  (\prod_j  q^{i_j-j}) . \]

\subsection{Finitization}
As in \ref{alg comp}, consider $ \sigma^N=0 $, and the differentials taken to be
\[\sigma,\sigma^{N-1}, \sigma,\sigma^{N-1}, \cdots, \]
with a $ \sigma $ passing the middle two columns.
When we set $ N\mapsto \infty$, we have the $ q $-dimension to be the $ q $-Catalan numbers. 

We consider the modified Euler characteristic
 \begin{equation} \label{qchi}
 \chi_q = \sum_{i\in \mathbb{Z}} (-1)^iq^ { f_{P(N)}(i) }  \qdim_i  \end{equation}
where $ f_{P(N)}(i) $ is an integer-valued function depending on our value of $ N $, and the partition into two parts of $ N $, $ P(N) $ that we consider. To obtain $ q $-Catalan numbers in the limit, note that  $ P(N) $ should be taken to be $ 1+(N-1) $.
A natual way to do so is to aim for $ q $-analogues of proposition \ref{prop1}.

\begin{prop} \label{prop3}
	For $N=3 $ and its unique partition $ 3=2+1 $, we have \[f_{1+2}(i)= \frac{3i^2-i}{2},\] so that the alternating sum produces $ 0, \pm1.$\\
	For  $ N=5 $, partition $ 5=1+4 $, \[f_{1+4} (i)= \frac{5i^2-3i}{2}, \] so that the alternating sum produces one version of the $ q $-fibonacci number with positive and negative signes, $ \pm F_n $. Here
	\[F_n= F_{n-1}+q^{n-1} F_{n-2} , \] where $ F_0=0, F_1=1; $\\
	
	for  $ N=5 $, partition $ 5=2+3$, \[f_{2+3} (i)= \frac{5i^2-i}{2}, \]
	so that the alternating sum produces another version of the $ q $-fibonacci number with positive and negative signes, $ \pm F'_n $. Here
	\[F'_n= F'_{n-1}+q^{n-2} F'_{n-2} , \] where $ F'_0=0, F'_1=1. $
\end{prop}

We note that for $ N=3 $, the correction terms are the pentagonal numbers (sequence A001318 in the OEIS); the two for $ N=5 $ are the exponents from two Roger-Ramanujan identities. 

\subsection{Categorification Conjectures}

Comparing proposition \ref{prop3} and \ref{prop1}, we see that for $ N=3,5 $, the deformation, with modifications on the $ q$-dimension via the coefficients of some sequences works well. The two different $ q$-Fibonacci numbers arise from the two distinct partitions of $ 5 $. In fact, works such as \cite{SA} show that there are more identities of Roger-Ramanujan type, and each of them gives a sequence of our $ q$-Euler character modifications for higher $ N $. 

This gives a possibe appoach to the categorification of both an entire set of $ q $-deformed sequences such as the ones appearing in proposition \ref{prop1} and the ones arise from the exponents of identities of the Roger-Ramanujan type.

\begin{conj}
The modified Euler characteristic \ref{qchi} with $  f_{P(N)} $ the exponents from identities of the Roger-Ramanujan type gives $ q $-deformed sequences arise from the usual Euler character as in proposition \ref{prop1}. 
\end{conj}

\begin{conj}
	There exists constructions of algebraic complexes such that their dimension matches the entries in the higher Pascal's triangles in the same way as for the Catanlan numbers. Moreover, ther exists $ q $-deformations of them natually give rise to $ q $-deformed sequences such as those in proposition \ref{prop2}.
\end{conj}

\end{document}